 \documentclass[final,1p]{elsarticle}
\usepackage{amssymb}
\usepackage{amsthm}
\usepackage{bm}
\usepackage{amsmath}
\usepackage{enumerate}
\usepackage{enumitem}
\usepackage[caption=false]{subfig}
\newtheorem{theorem}{Theorem}[section]

\journal{Wave Motion}

\begin{document}

\begin{frontmatter}



\title{Numerical Computation of  Solitary Wave Solutions of the Rosenau Equation}


\author[label1]{H. A. Erbay\corref{cor1}}
    \ead{husnuata.erbay@ozyegin.edu.tr}
    \address[label1]{Department of Natural and Mathematical Sciences, Faculty of Engineering, Ozyegin University,  Cekmekoy 34794, Istanbul, Turkey}
    \cortext[cor1]{Corresponding author}

\author[label1]{S. Erbay}
    \ead{saadet.erbay@ozyegin.edu.tr}

\author[label2]{A. Erkip}
    \ead{albert@sabanciuniv.edu}
    \address[label2]{Faculty of Engineering and Natural Sciences, Sabanci University, Tuzla 34956, Istanbul,  Turkey}

\begin{abstract}
    We construct numerically solitary wave solutions of the Rosenau equation using the Petviashvili iteration method.  We first summarize the theoretical results available in the literature for the existence of solitary wave solutions. We then apply two numerical algorithms based on the Petviashvili method for solving the Rosenau equation with single or double power law nonlinearity.  Numerical calculations   rely on a uniform discretization of a finite computational domain.  Through some numerical experiments we observe that the algorithm converges rapidly and it is robust to  very general forms of the initial guess.
\end{abstract}

\begin{keyword}
    Rosenau equation \sep Petviashvili iteration method \sep   Benjamin-Bona-Mahony equation \sep solitary waves

    \MSC[2010] 35Q53 \sep   65M99 \sep 74J35 \sep 74S30
\end{keyword}

\end{frontmatter}


\setcounter{equation}{0}
\section{Introduction}\label{sec:sec1}

This paper is concerned with the construction of the solitary wave solutions numerically for the Rosenau equation
\begin{equation}
    u_{t}+u_{x}+u_{xxxxt}+(g(u))_x=0,  \label{eq:re}
\end{equation}
where $u=u(x,t)$ is a real valued function and  $g$ is a smooth nonlinear function to be specified later. The Rosenau equation  models, within a quasi-continuum framework, a unidirectional propagation of  longitudinal  waves on the one-dimensional  crystal lattice which is assumed to be dense  \cite{Rosenau1988}. It is well-known that nonlinear dispersive wave equations may possess smooth, symmetric, localized traveling-wave solutions with monotonic tails (called solitary waves). An explicit solitary wave solution of the form
\begin{equation}
  u(x,t)=\mbox{sech}(x-{1\over 2}t)~ \label{eq:sech}
\end{equation}
was presented in \cite{Park1992}  for  the particular form
\begin{equation}
   u_{t}+ u_x+u_{xxxxt}-30 u^2 u_x+60 u^4 u_x=0,  \label{eq:re-q}
\end{equation}
which is a special case of the Rosenau equation (\ref{eq:re}) with the quintic polynomial  $~g(u)=-10u^3+12u^5$. To the best of the knowledge of the authors, (\ref{eq:sech}) is the only explicit solitary wave solution of the Rosenau equation, given in the literature. In \cite{Esfahani2014}, an explicit solitary wave solution to the Rosenau-Regularized Long Wave  (Rosenau-RLW) equation with quadratic nonlinearity
\begin{equation}
  u_{t}+ u_x+\alpha u_{xxx}+ u_{xxxxt}+uu_x=0  \label{eq:rre}
\end{equation}
was given in the form
\begin{equation}
  u(x,t)={35\over 12}(c-1)\mbox{sech}^4\left ( \sqrt{{13(1-c)\over 144 c\alpha}}(x-ct) \right ), \label{eq:rre-sw}
\end{equation}
with $c=169/(169-36 \alpha^2)>1$ and $\alpha<0$.   Note that, when $\alpha=0$, (\ref{eq:rre}) reduces to the Rosenau equation with $g(u)=u^2/2$.  However, the  solitary wave solution   (\ref{eq:rre-sw}) is not defined in the limit $\alpha \to 0^{-}$.

The existence of  solitary wave solutions to the Rosenau equation  has been proved in \cite{Zeng2003} (see Section \ref{sec:sec2} for more details). However, it is still an open question whether explicit  form of the solitary wave solution of the Rosenau equation with a single power type nonlinear term  $g(u)=u^{p+1}/(p+1)$ ($p\ge 1$ is an integer) can be found.  So presently, it is not possible to have an idea both about the solution profiles for various power nonlinearities and about  how the amplitude changes with the degree of nonlinearity. Furthermore, our expectation is that the existence of the higher-order dispersive term in (\ref{eq:re}) may cause solitary waves with non-monotonic tails. Naturally these considerations motivate us to provide a numerical scheme to construct numerically the solitary wave solutions of the Rosenau equation with power nonlinearities. For this aim  we use  the Petviashvili iteration method \cite{Petviashvili1976}, which has been widely used to compute solitary wave solutions of nonlinear dispersive wave equations (see for instance \cite{Ablowitz2005,Fibich2006,Lakoba2007,Duran2014, Borluk2017,Duran2018,Dougalis2019} and references therein). The main idea behind this method is to introduce a stabilizing factor to the fixed-point iteration scheme used to get solitary wave solutions. The main advantages of the Petviashvili iteration method are   fast convergence in most of the cases and   robustness to the choice of the initial guess. We refer the reader to \cite{Pelinovsky2004,Olson2016} for convergence analysis of the method.

Our numerical experiments have led us to the conclusion that the existence of higher-order dispersive term in (\ref{eq:re}) leads to the existence of the localized waves with   non-monotonic tails.  For clarity, we underline that our discussion and reported results in the following sections do not cover the explicit solution (\ref{eq:sech}). The reason is that our numerical approach is based on wave speeds greater than unity whereas the wave speed $1/2$ of (\ref{eq:sech}) is less than unity.

The paper is structured as follows. In Section \ref{sec:sec2}, we present  some theoretical results involving both local well-posedness of the Cauchy problem  and the existence of solitary wave solutions to (\ref{eq:re}). In Section \ref{sec:sec3}, we apply the Petviashvili iteration method to  compute   solitary wave solutions of (\ref{eq:re}). Finally in Section 4 we perform some numerical experiments to test the performance of the proposed algorithm.


\setcounter{equation}{0}
\section{Preliminaries}\label{sec:sec2}

\subsection{Local Well-Posedness of the Cauchy Problem}
In \cite{Albert1991} the  well-posedness of the initial-value problem for the following  class of nonlinear dispersive equations
\begin{equation}
     u_t+u_x + u^p u_x+Mu_t=0 \label{eq:zeng}
\end{equation}
has been considered. Here    $M$ is a Fourier multiplier operator in the $x$ variable defined by
\begin{equation}
    \widehat{M v}(\xi) = m(\xi) \widehat{v}(\xi),
\end{equation}
where the symbol $\widehat{~}$ denotes the Fourier transform. When $p=1$ and $M=-D_x^2$  for which $m(\xi)=\xi^2$ (where $D_{x}$ represents the partial derivative with respect to $x$), (\ref{eq:zeng}) becomes  the Benjamin-Bona-Mahony (BBM) equation \cite{Benjamin1972},
\begin{equation}
    u_{t}+u_{x}-u_{xxt}+u^{p}u_{x}=0,  \label{eq:bbm}
\end{equation}
that has been widely used  to model unidirectional surface water waves with small amplitudes and long wavelength.  We note that the Rosenau equation (\ref{eq:re}) with $g(u)=u^{p+1}/(p+1)$ is a member of the class (\ref{eq:zeng}) when $M=D_x^4$.  Hence the following well-posedness result on the Sobolev space $H^{s}(\mathbb{R})$ given in  \cite{Albert1991} for (\ref{eq:zeng}) is also valid for  the present form of (\ref{eq:re}):
\begin{theorem}[\cite{Albert1991}]  \label{theo2}
    Consider the initial value problem for (\ref{eq:zeng}) with $u(x,0)=u_0(x)$. Suppose that  the symbol $m$ of $M$ satisfies
    \begin{eqnarray*}
        && c_1 \vert \xi\vert^\mu \leq m(\xi)\quad \mbox{for all~} \vert\xi\vert \geq 1\\
        && m(\xi)\leq c_2(1+\vert \xi\vert)^\nu \quad \mbox{for all~} \xi\in \mathbb{R},
    \end{eqnarray*}
    where $c_1$ and $c_2$ are positive constants and $1\leq \mu\leq \nu$.  If $u_0\in H^s(\mathbb{R})$ where $s > 1/2$,  then there exists a $T_0>0$ depending only upon $\Vert u_0\Vert_{H^s}$ such that the initial value problem admits a unique solution which for any $T < T_0$, lies in $C^{\infty}([0, T];H^s(\mathbb{R}))$.  Moreover, the correspondence that associates to initial data  $u_0$   the unique solution $u$ is continuous from $H^s(\mathbb{R})$ to $C^k([0, T];H^s(\mathbb{R}))$ for any $T < T_0$ and any finite value of $k$. The existence time $T_0$ depends inversely on
 $\| u_0\|_{H^s}$  and $T_0\to \infty$ as $\vert| u_0\vert|_{H^s}\to 0$.
\end{theorem}

In \cite{Erbay2019}, the  well-posedness of the initial value problem for the following  class of nonlocal nonlinear dispersive equations
\begin{equation}
    u_{t} +\left(\beta \ast  f(u)\right)_{x}=0  \label{eq:nl-eq}
\end{equation}
has been established.  In (\ref{eq:nl-eq}), $f$ is a sufficiently smooth function with $f(0)=0$, $\beta$ is a general kernel function and  $\beta \ast v$ denotes the convolution of $\beta$ with  the function $v$
\begin{displaymath}
    (\beta \ast v)(x)= \int_{\mathbb{R}} \beta(x-y)v(y)\mbox{d}y.
\end{displaymath}
The above class of convolution-type equations involves many well-known wave equations as particular cases. For instance, in the limit case in which the kernel function $\beta$ is taken as the Dirac delta function $\delta(x)$, (\ref{eq:nl-eq}) reduces to the hyperbolic conservation equation $u_{t} +(f(u))_{x}=0$. Similarly, if  $\beta(x)=e^{-|x|}/2$  and if $f(u)=u+u^{p+1}/(p+1)$ with  $p\geq 1$,  (\ref{eq:nl-eq}) reduces to  the BBM equation (\ref{eq:bbm}). Note that the exponential kernel $\beta(x)=e^{-|x|}/2$ is the Green's function for the differential operator $1-D^{2}_{x}$. On the other hand, if the kernel function $\beta$ is chosen as the  Green's function of the differential operator $1+D^{4}_{x}$,
\begin{equation}
    \beta(x)= {1\over {2\sqrt 2}}e^{-{\vert x\vert\over \sqrt 2}}\Big( \cos\big({{\vert x\vert}\over {\sqrt 2}}\big) + \sin \big({{\vert x\vert}\over {\sqrt 2}}\big) \Big ), \label{eq:ros-ker}
\end{equation}
and if $f(u)=u+g(u)$,  (\ref{eq:nl-eq}) reduces to  the  Rosenau equation (\ref{eq:re}).

Assuming that $\beta$ is integrable or more generally a finite measure on $\mathbb{R}$ with positive Fourier transform, the convolution with the kernel function $\beta$ is a positive bounded operator on the Sobolev space $H^{s}(\mathbb{R})$.  The well-posedness of the Cauchy problem for  (\ref{eq:nl-eq}) given in \cite{Erbay2019} is also valid for the Rosenau equation:
\begin{theorem}[\cite{Erbay2019}] \label{theo3}
    Consider the initial-value problem for (\ref{eq:nl-eq}) with $u(x,0)=u_0(x)$.  Let $s > 1/2$, $f\in C^{[s]+1} (\mathbb{R})$ with $f(0)=0$, $\beta\in L^1(\mathbb{R})$ and $ \beta'=\mu$ is a finite measure on $\mathbb{R}$. For a given  $u_0\in H^s(\mathbb{R})$, there is some $T > 0$ so that the initial-value problem   is locally well-posed with solution $u \in C^1\left( [0,T];H^s(\mathbb{R})\right ).$
\end{theorem}

\subsection{Existence and Stability of Solitary Wave Solutions}

The existence and stability of the solitary wave solutions of (\ref{eq:zeng}) have been investigated in \cite{Zeng2003}.     Using the Concentration-Compactness Method  \cite{Lions1984}, the author has established the existence and stability of solitary wave solutions of (\ref{eq:zeng}) provided that $p$ and $m(\xi)$ satisfy the following conditions:
\begin{eqnarray*}
    & A1.&  \text{There exist positive constants} ~c_1 ~\text{and}~ r>{p\over 2} ~\text{such that}~
        m(\xi)\leq c_1 |\xi|^r~\text{for}~ |\xi|\leq 1;\\
    & A2.&  \text{There exist positive constants} ~c_2 ~\text{and}~ c_3 ~\text{ and}~s\geq 1~\text{such that}~\\
        &&        c_2|\xi|^s\leq m(\xi)\leq c_3 |\xi|^s~\text{for}~ |\xi|\geq 1;\\
    & A3.& m(\xi)\geq 0 ~\text{for all values of }~\xi;\\
    & A4.& m(\xi)~\text{is four times differentiable for all non-zero values of }~\xi, ~\text{and for each}~\\
        &&  j\in \{0,1,2,3,4\}~\text{there exist positive constants} ~c_4 ~\text{and}~ c_5 ~\text{such that}
\end{eqnarray*}
\begin{eqnarray*}
    && \left \vert  {d^j\over d \xi^j}   \left({m(\xi)-m(0)\over \xi}\right )\right \vert
        \leq {c_4\over \vert \xi\vert^j}~~\text{for }~~0<\vert\xi\vert \leq 1,\\
    && \left \vert  {d^j\over d \xi^j}   \left({\sqrt{m(\xi)}\over \xi^{s\over 2}}\right )\right \vert
        \leq {c_5\over \vert\xi\vert^j}~~\text{for }~~\vert\xi\vert\geq 1.
\end{eqnarray*}
Recall that the Rosenau equation  (\ref{eq:re}) is a member of the above class (\ref{eq:zeng}) with  $M=D_x^4$ (that is, $m(\xi)=\xi^4$). It is easy to check that the assumptions $A1-A4$ are satisfied when $p<8$ and $s=4$.  Thus we conclude that the Rosenau equation  has solitary waves solutions when the nonlinear term is of the form $g(u)=u^{p+1}/(p+1),~p\ge 1$. On the other hand, except the solution (\ref{eq:sech}) given in \cite{Park1992} for a very special  nonlinearity, no explicit solution of the Rosenau equation is known. The above results on the existence of the solitary waves solutions to  (\ref{eq:zeng}), in particular the Rosenau equation, motivate us to construct solitary waves solutions numerically.

\setcounter{equation}{0}
\section{The Petviashvili Method}\label{sec:sec3}

In this section we briefly discuss the Petviashvili method and apply the method to construct  solitary wave solutions of (\ref{eq:re}). The Petviashvili iteration method was introduced in \cite{Petviashvili1976} to construct numerically solitary wave solutions of one or higher dimensional nonlinear wave equations. Recently  a number of articles have reported that the method  can be successfully applied to nonlinear dispersive wave equations  \cite{Ablowitz2005,Fibich2006,Lakoba2007,Duran2014, Borluk2017,Duran2018,Dougalis2019}.  In  \cite{Pelinovsky2004},  the convergence of the iteration scheme was proved  for a general nonlinear wave equation with single power nonlinearity and    the conditions necessary to achieve the optimal convergence rate were found. The robustness of the method to the initial guess  has been discussed in \cite{Olson2016}.

For the Rosenau equation (\ref{eq:re}) with  power nonlinearity  $\displaystyle g(u)=u^{p+1}/(p+1)$, let us consider traveling wave solutions $u(x,t)=\phi(x-ct)$ where  $c$ is the constant wave speed. Assuming that $\phi(x)$ and its derivatives tend to zero as  $\vert x\vert\to \infty$, after one integration, substitution of the traveling wave solution  into (\ref{eq:re}) leads to the equation
\begin{equation}
    \mathcal{L}\phi= (cD^4+(c-1))\phi={\phi^{p+1} \over {p+1}},  \label{eq:sw}
\end{equation}
where $D$ and $\mathcal{L}$ denote differentiation and the fourth-order differential operator, respectively.  For $c>1$ the operator $\mathcal{L}$ is invertible, and we may write (\ref{eq:sw}) as
\begin{equation}
    \phi= {1\over p+1}\mathcal{L}^{-1} \phi^{p+1} \label{eq:sw1}.
\end{equation}
In the Fourier space, (\ref{eq:sw1}) becomes
\begin{equation}
    \widehat{\phi}(\xi)= {1\over p+1}{\widehat{\phi^{p+1}}(\xi) \over {l(\xi)}}, ~~~~l(\xi)=c \xi^4+c-1, \label{eq:sw2}
\end{equation}
 where $l(\xi)$ is the symbol of $ \mathcal{L}$ and  the   Fourier transform and its inverse are defined as
\begin{equation}
    \widehat{\phi}(\xi)= \int_\mathbb{R}   \phi(x) e^{-i\xi x} dx,\quad  \quad
    \phi(x)={1\over 2 \pi}\int_\mathbb{R} \widehat{\phi}(\xi) e^{i\xi x} d\xi.
\end{equation}
It is worth pointing out that under the assumption  $c>1$ we have $l(\xi)>0$  for any $\xi \in \mathbb{R}$. Multiplying (\ref{eq:sw2}) by the complex conjugate $(\widehat{\phi})^{*}$  and integrating over $\mathbb{R}$ gives
\begin{equation}
    \big\langle l(\xi)\widehat{\phi}(\xi),~  (\widehat{\phi})^{*}(\xi)\big\rangle
        = {1\over {p+1}}\big\langle \widehat{\phi^{p+1}}(\xi),~ (\widehat{\phi})^{*}(\xi) \big\rangle \label{eq:sw3},
\end{equation}
where the symbol  $\langle~\cdot~, ~\cdot~\rangle$ is used to denote the standard inner product in Fourier space, defined by
\begin{displaymath}
  \big\langle f(\xi), g(\xi) \big\rangle = \int_\mathbb{R} f(\xi) g^{*}(\xi) d\xi .
\end{displaymath}

The Petviashvili method suggests that (\ref{eq:sw2}) can be solved by considering the following iteration scheme
\begin{equation}
    \widehat{\psi}_{n+1}(\xi)={1\over {p+1}} {\widehat{\psi_n^{p+1}}(\xi)\over {l(\xi)}},~~~~ \quad n=0,1,\cdots,
\end{equation}
where $\psi_n(x)$ represents the approximation at the $n$th iteration to $\phi(x)$. Even if a non-trivial fixed point $\widehat{\phi}(\xi)$ exists,  this iteration scheme may diverge to trivial fixed points $\phi=\infty$ or  $\phi=0$. To resolve this problem, a  modified iteration scheme involving a stabilizing  factor was introduced in \cite{Petviashvili1976} by Petviashvili. The main idea of the modified iteration scheme  is to renormalize the approximate solution at each iteration step so that the identity (\ref{eq:sw3}) is satisfied. We define the stabilizing factor $P_{n}$ as
 \begin{equation}
    P_n=(p+1){\big\langle l(\xi)\widehat{\psi}_n(\xi),~(\widehat{\psi}_n)^{*}(\xi) \big\rangle \over
        \big\langle  \widehat{\psi_n^{p+1}}(\xi),~ (\widehat{\psi_n})^{*}(\xi)\big\rangle},~~~~ \quad n=0,1,\cdots \label{eq:Pstab}
\end{equation}
As expected, the $n$th iterate $ \widehat{\psi}_n$  does not satisfy the identity (\ref{eq:sw3}). However, we note that  $P_{n} \rightarrow 1$ as $n \rightarrow \infty$ when the scheme converges, that is, when $\psi_{n} \rightarrow \phi$ as  $n \rightarrow \infty$.
Following Petviashvili \cite{Petviashvili1976} we introduce the stabilizing factor $P_n$ into the iteration as follows
\begin{equation}
    \widehat{\psi}_{n+1}(\xi)= {(P_n)^{\theta}\over {p+1}} {\widehat{\psi_n^{p+1}}(\xi) \over {l(\xi)}},   \quad n=0,1,\cdots. \label{eq:meth2}
\end{equation}
As it was pointed out by Petviashvili,  the fastest convergence of iterations occurs if the parameter $\theta$ is taken as $\theta=(p+1)/p$. The origin of this optimal value of $\theta$ was shown  in \cite{Pelinovsky2004}. With a similar argument we use this optimal value of $\theta$ in our numerical experiments below.

Since the trivial (zero) solution is a solution for (\ref{eq:sw}), the iteration scheme will not converge to a nonzero solution if the starting function $\psi_{0}(x)$ lies in the domain of attraction of the zero solution. Additionally,  the scheme  may not converge at all if the starting function $\psi_{0}(x)$ is not chosen sufficiently close to the solitary wave solution.  Through the numerical experiments in the next section, for various forms of the starting function,  we observe  that the sequence $\widehat{\psi}_n(\xi)$ converges rapidly to the fixed point of (\ref{eq:meth2}) and consequently to the solitary wave solution of the Rosenau equation (\ref{eq:re}).

We remark that the above iteration scheme is valid for  a homogeneous power nonlinearity of degree ($p+1$). For more complicated nonlinearities we may need more than one stabilizing factor. In the next section, we also discuss one such example, the cubic-quintic Rosenau equation.

\setcounter{equation}{0}
\section{Numerical Experiments}\label{sec:sec4}

In this section we perform  some numerical experiments to show the convergence of  the Petviashvili method  and to compare the numerical solutions obtained   with the exact solutions available. All the numerical simulations   are carried out using Matlab. In our setup, we replace the infinite domain by a finite computational domain  $[-L,L]$ and  discretize it  using $2N$ equally-spaced  subintervals with the grid spacing $h=L/N$. The discrete Fourier transform and its inverse are performed using the so-called fast Fourier transform (via the  Matlab routines \verb"fft" and \verb"ifft", respectively). The integrals are computed via  trapezoidal integration (the  Matlab routine \verb"trapz").

There are different sources of  error in our computations: the domain truncation error resulting from the restriction of the infinite interval to a finite interval of computation, the discretization error resulting from the consideration of a finite number of grid points and the algorithmic error resulting from the consideration of a finite number of fixed point iteration. For sufficiently large computational interval, the domain truncation error    has no significant effect on the  numerical results, since the solitary wave  solutions to be obtained will decay rapidly to zero for $\vert x\vert \rightarrow \infty$.  Furthermore, the number of grid points is chosen large enough to minimize the discretization errors.  So the focus of the numerical experiments will be the algorithmic error rather than  the domain truncation errors or the discretization errors. Following the literature we may measure the algorithmic errors in three different ways.  The   residual error $E_{n}^{(r)}$,  the iteration error $E_{n}^{(s)}$ based on the difference from unity of the stabilizing factor and the iteration error $E_{n}^{(a)}$ based on the difference between two successive  approximations are calculated, respectively, as
\begin{equation}
    E_{n}^{(r)}= \big \Vert \mathcal{L}\psi_{n}-{\psi_{n}^{p+1}\over {p+1}}\big \Vert_{L^{\infty}}, ~~~~
    E_{n}^{(s)}= \big \vert 1-P_{n}\big \vert, ~~~~
    E_{n}^{(a)}= \big \Vert \psi_{n+1}-\psi_{n}\big \Vert_{L^{\infty}}. \label{eq:error}
\end{equation}
However, since, in the numerical experiments presented below, we observed that the main characteristics of the results for $E_{n}^{(s)}$ and $E_{n}^{(a)}$ are almost identical, we will present only the results for $E_{n}^{(s)}$ below for simplicity. In all the experiments we set the tolerance for both the residual error $E_{n}^{(r)}$ and the iteration error $E_{n}^{(s)}$ to be less than $10^{-14}$. Furthermore we set the maximum number of iterations to be performed in the Petviashvili algorithm to 500, but in most cases  we only need 40 to 50 iterations to obtain the desired accuracy.

\subsection{The  BBM Equation}

To be able to present the rapid convergence of the Petviashvili iteration scheme,  we first apply the  method to the generalized BBM equation (\ref{eq:bbm}) for which an explicit solitary wave solution is available. The solitary wave solution of the generalized BBM equation (\ref{eq:bbm}) which is a member of the class (\ref{eq:nl-eq}) with the exponential kernel is given by.
\begin{equation}
    u(x,t)=A~\Big(\text{sech}^{2}\big(B(x-ct-x_{0}) \big)\Big)^{1/p}, \label{eq:solitary}
\end{equation}
with $A=\big((p+1)(p+2)(c-1)/2\big)^{1/p}$, $B=(p/2)(1-1/c)^{1/2}$ and $c >1$ \cite{Bona2000}.  The solitary wave  (\ref{eq:solitary}) is initially located at $x_{0}$ and propagates to the right  with the constant wave speed $c$. From now on we fix $c=1.8$.

\begin{figure}[h!]
    \centering
    \subfloat[Exact and numerical profiles]{\label{fig:1a}
        \includegraphics[height=0.24\textheight,scale=1.50,keepaspectratio]{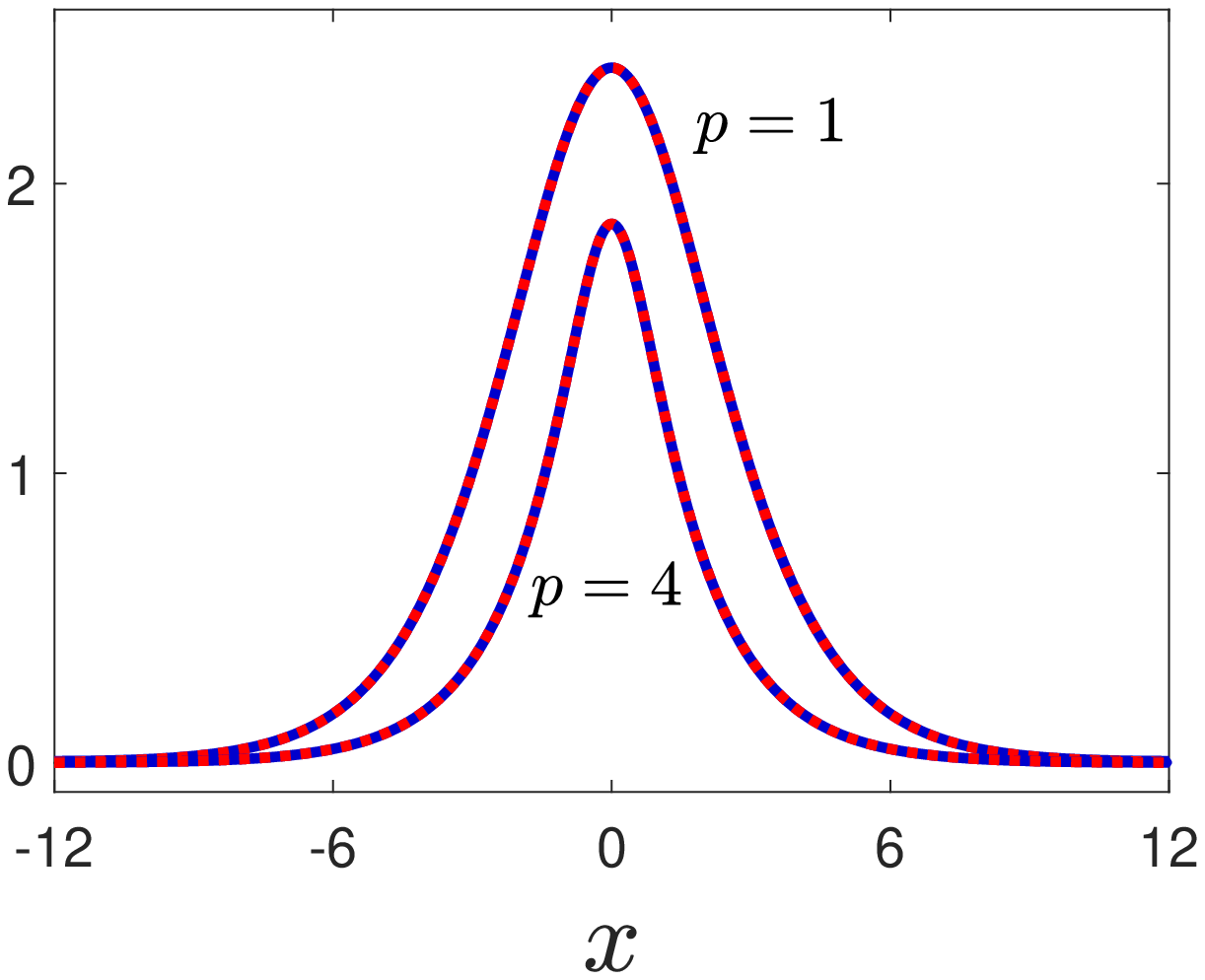}}\hspace*{10pt}
    \subfloat[Residual and iteration errors]{\label{fig:1b}
       \includegraphics[height=0.24\textheight,scale=1.5,keepaspectratio]{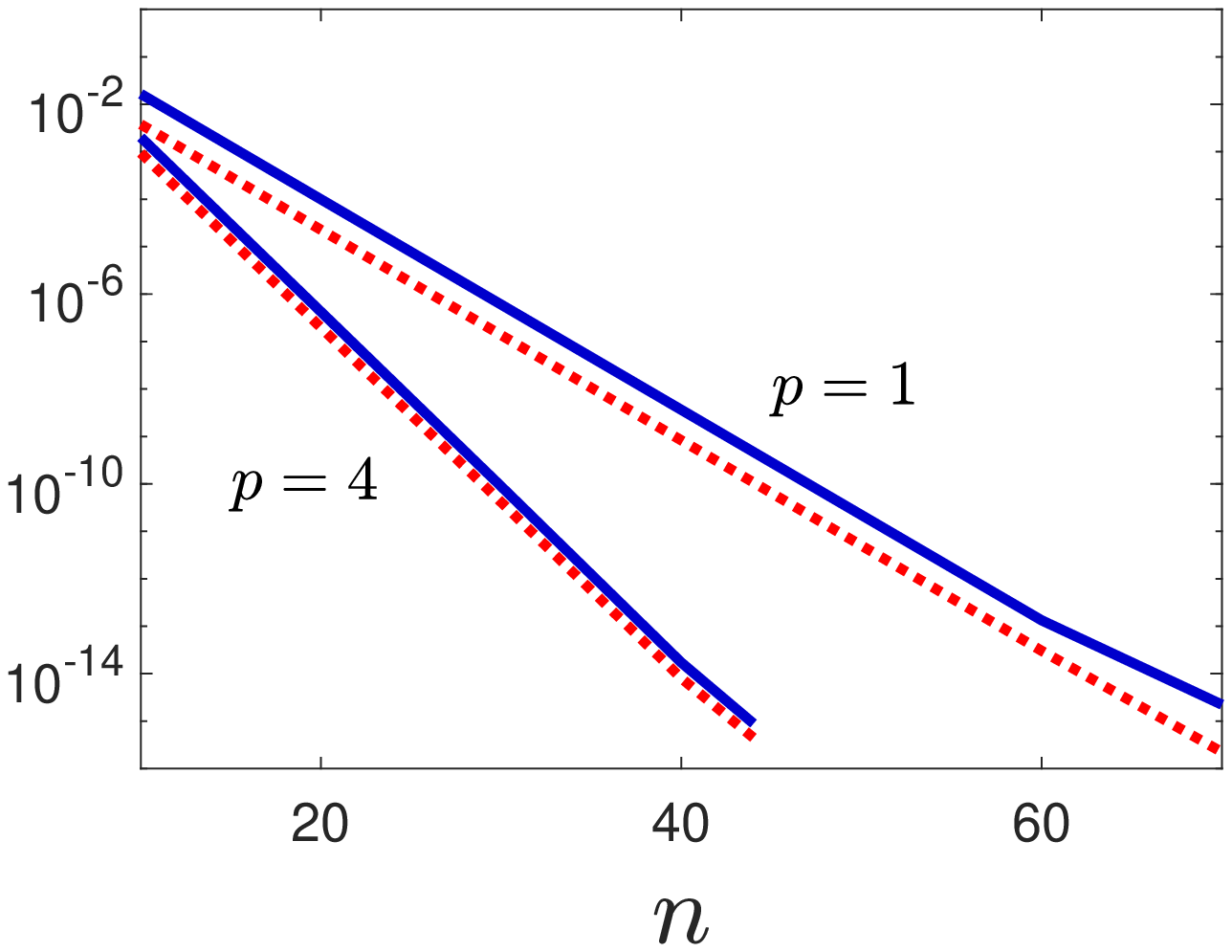}}
    \caption{(a) The exact and the numerical solution profiles  of the BBM equation (quadratic nonlinearity ($p=1$) and quintic nonlinearity ($p=4$)) for the Gaussian starting function. The dotted (red) line and the solid (blue) line represent the exact and numerical solutions, respectively. The curves corresponding to the exact and numerical solutions coincide almost completely.   (b) Variation of  the iteration  error $E_{n}^{(s)}$ and the residual error $E_{n}^{(r)}$ with the iteration number $n$ for the BBM equation (quadratic nonlinearity ($p=1$) and quintic nonlinearity ($p=4$)) for the Gaussian starting function. The dotted (red) line and the solid (blue) line represent $E_{n}^{(s)}$ and $E_{n}^{(r)}$, respectively. (The wave speed $c=1.8$, the computational domain $[-12, \, 12]$ and the mesh size   $0.05$ are used in all computations.)}
    \label{fig:Fig1}
\end{figure}

We remark that the Petviashvili algorithm defined by (\ref{eq:Pstab}) and (\ref{eq:meth2}) is also valid for  the BBM equation if we take $l(\xi)=c \xi^2+c-1 $. We take the initial guess $\psi_{0}$  as the Gaussian function $\displaystyle    \psi_{0}(x)=e^{-x^{2}}$. In all the numerical experiments related to the BBM equation we fix the size of the computational domain and the mesh size as $-12\leq x \leq 12$ and $h=0.05$ (which corresponds to 480 discrete Fourier modes), respectively.  We now apply the Petviashvili method for the parameter values $p=1$ and $p=4$ and plot the exact and numerical solutions in Figure \ref{fig:1a}. For both values of $p$, we observe that the numerical solution has the same amplitude and waveform as that of the exact solution and overall there is a very good agreement between the two curves. In order to verify the convergence of the method we now conduct a sequence of numerical experiments for different values of iteration number $n$. Using a semi-logarithmic scale, in Figure \ref{fig:1b} we  present variation of the residual error $E_{n}^{(r)}$ and the iteration error $E_{n}^{(s)}$ with the number of iterations ($n$) in the fixed point algorithm. We observe that the iteration scheme rapidly converges to the solitary wave solution. We make similar calculations for various values of $p$ but we always get a complete agreement between the exact and numerical solutions.

\begin{figure}[h!]
    \centering
    \subfloat[Numerical profiles]{\label{fig:2a}
    \includegraphics[height=0.24\textheight,scale=1.50,keepaspectratio]{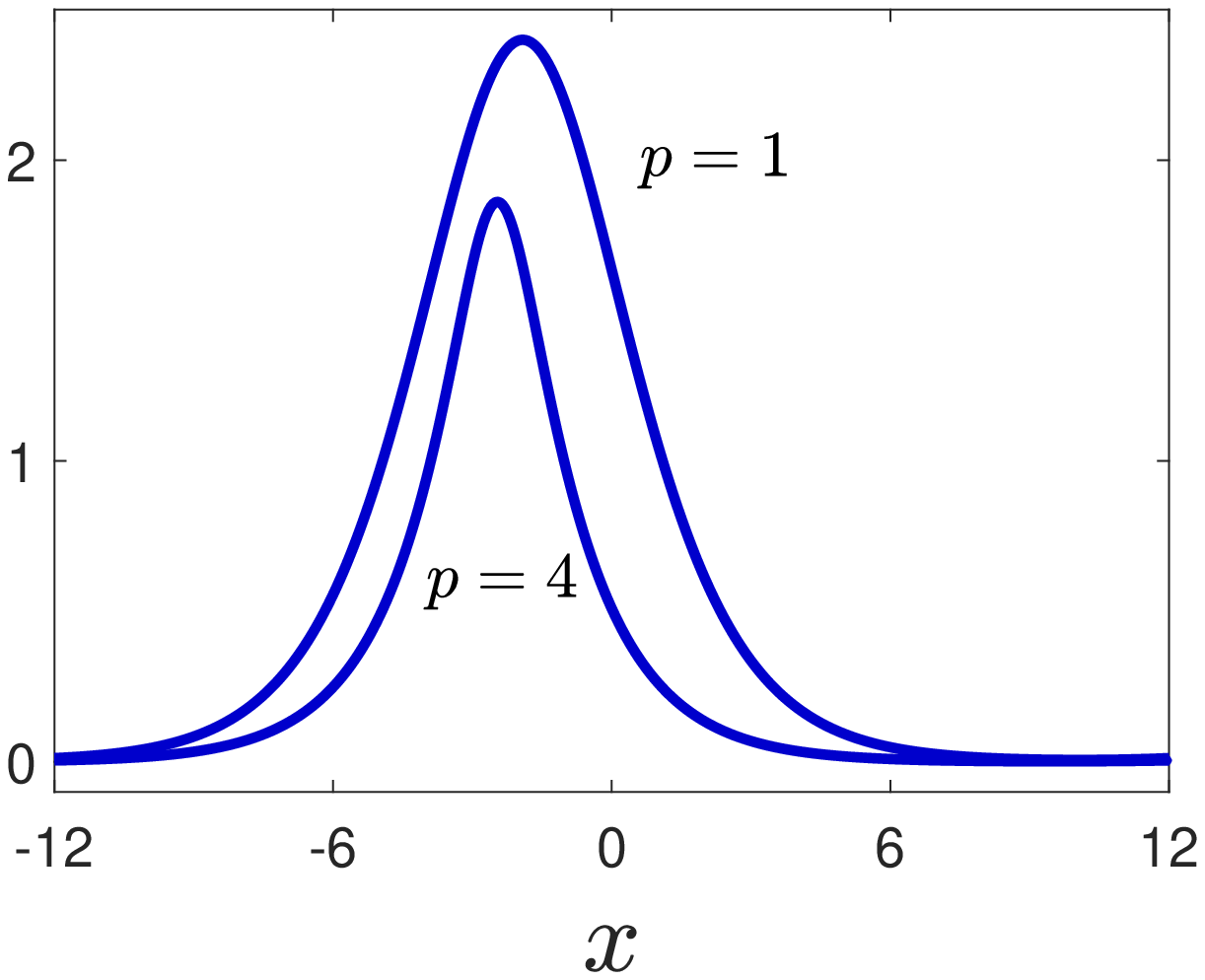}}\hspace*{10pt}
    \subfloat[Residual and iteration errors]{\label{fig:2b}
    \includegraphics[height=0.24\textheight,scale=1.50,keepaspectratio]{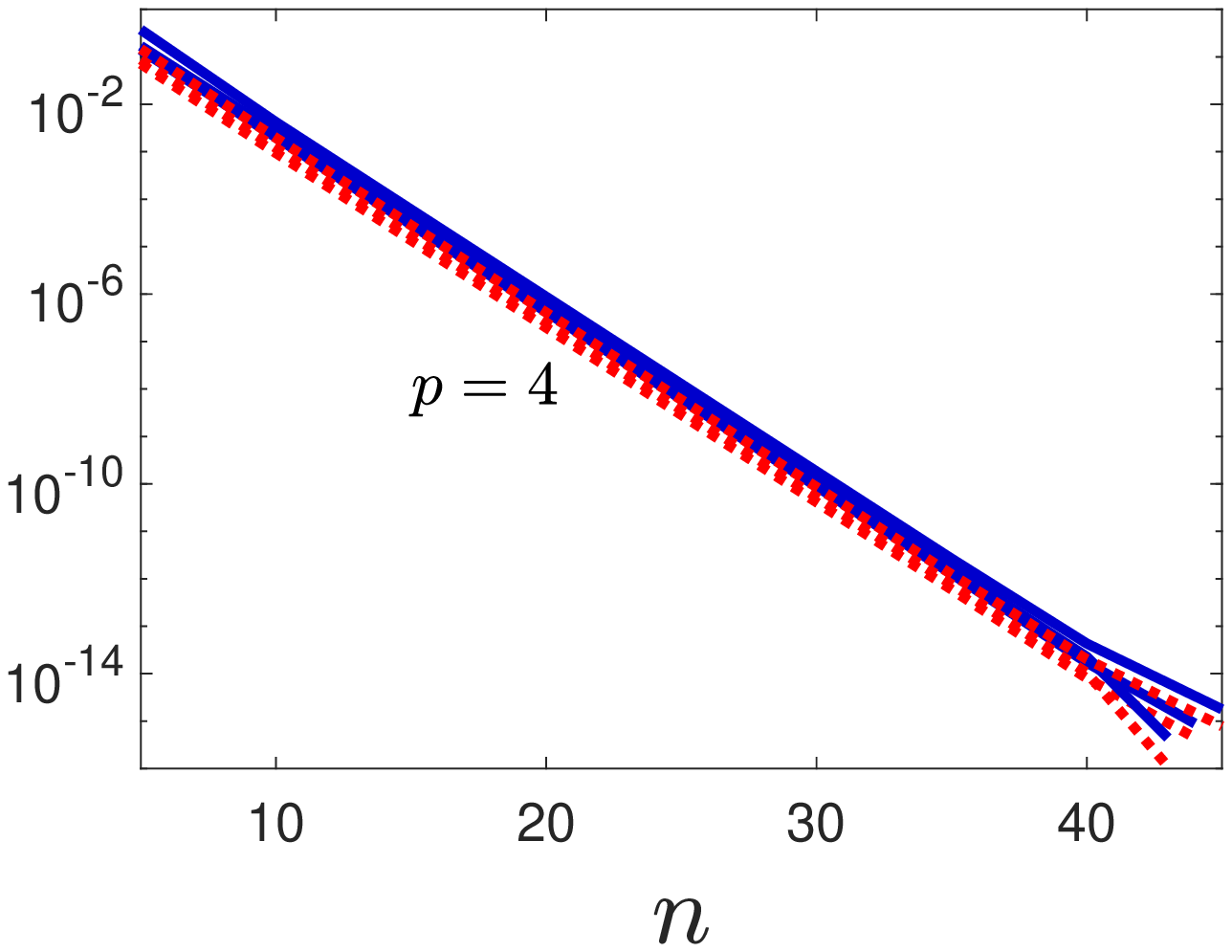}}
    \caption{(a) The numerical solution profiles of the BBM equation (quadratic nonlinearity ($p=1$) and quintic nonlinearity ($p=4$)) for the triangular starting function.     (b) Variation of  the iteration  error $E_{n}^{(s)}$ and the residual error $E_{n}^{(r)}$ with the iteration number $n$ for  three different types of the starting function: Gaussian function, nonsmooth exponential function and triangular function. The BBM equation is considered with  quintic nonlinearity ($p=4$). The dotted (red) line and the solid (blue) line represent $E_{n}^{(s)}$ and $E_{n}^{(r)}$, respectively. (The wave speed $c=1.8$, the computational domain $[-12, \, 12]$ and the mesh size   $0.05$ are used in all computations.)}
    \label{fig:Fig2}
\end{figure}

We now want to show that the  above numerical results are not very sensitive to the choice of the starting function $\psi_{0}(x)$.
In  the above numerical experiments we have taken the initial guess $\psi_{0}(x)$  as  the Gaussian function which is smooth and symmetric function about $x=0$. For both odd and even values of $p$  we have observed that  the algorithm converges rapidly to the solitary wave solution which  is symmetric  about the origin. In order to show that the Petviashvili method is extremely robust to the choice of the initial guess $\psi_{0}(x)$ we now repeat all the experiments for the nonsmooth symmetric function  $\displaystyle  \psi_{0}(x)=e^{-\vert x \vert}$. In those experiments we obtain almost the same numerical results with the same nature and we do not show the corresponding figures here. As the initial guess we now consider the asymmetric triangular function $\psi_{0}(x)$ defined by $\displaystyle  \psi_{0}(x)=x/3+2$ for $-6 \leq x \leq -3$, $\displaystyle  \psi_{0}(x)=-x/9+2/3$ for $-3 < x \leq 6$ and  $\displaystyle  \psi_{0}(x)=0$ otherwise   and repeat the experiments for various values of $p$. Again, we observe that the iteration scheme  converges very rapidly. However, this time, the iteration scheme  converges to a solitary wave solution which is a shift of the one in the previous case.  Since (\ref{eq:sw}) is translation invariant, we remark that if $\phi(x)$ is a solution, then the shifted function $\phi(x+a)$ is also a solution for any shift $a$.  We present the solution profiles for $p=1$ and $p=4$ in Figure \ref{fig:2a}.

We now want to evaluate the robustness of the iteration scheme to change in the starting function $\psi_{0}(x)$ in terms of errors. For this aim we consider the BBM equation  with  quintic nonlinearity ($p=4$) and then  compare the variation of the errors with the number of iterations for the above-mentioned three different initial guesses: Gaussian function, symmetric nonsmooth exponential function and asymmetric triangular function. We present the results for  the residual error $E^{(r)}_{n}$ and the iteration error $E^{(s)}_{n}$ with the number of iterations ($n$)  in Figure \ref{fig:2b}. We observe that the curves corresponding to three different forms of the starting function $\psi_{0}(x)$ are almost  indistinguishable from each other. We conclude that these experiments illustrate the robustness of the Petviashvili method.

\subsection{The  Rosenau Equation With Single Power Nonlinearity}
\begin{figure}[h!]
    \centering
    \subfloat[Numerical profiles]{\label{fig:3a} 
    \includegraphics[height=0.24\textheight,scale=1.50,keepaspectratio]{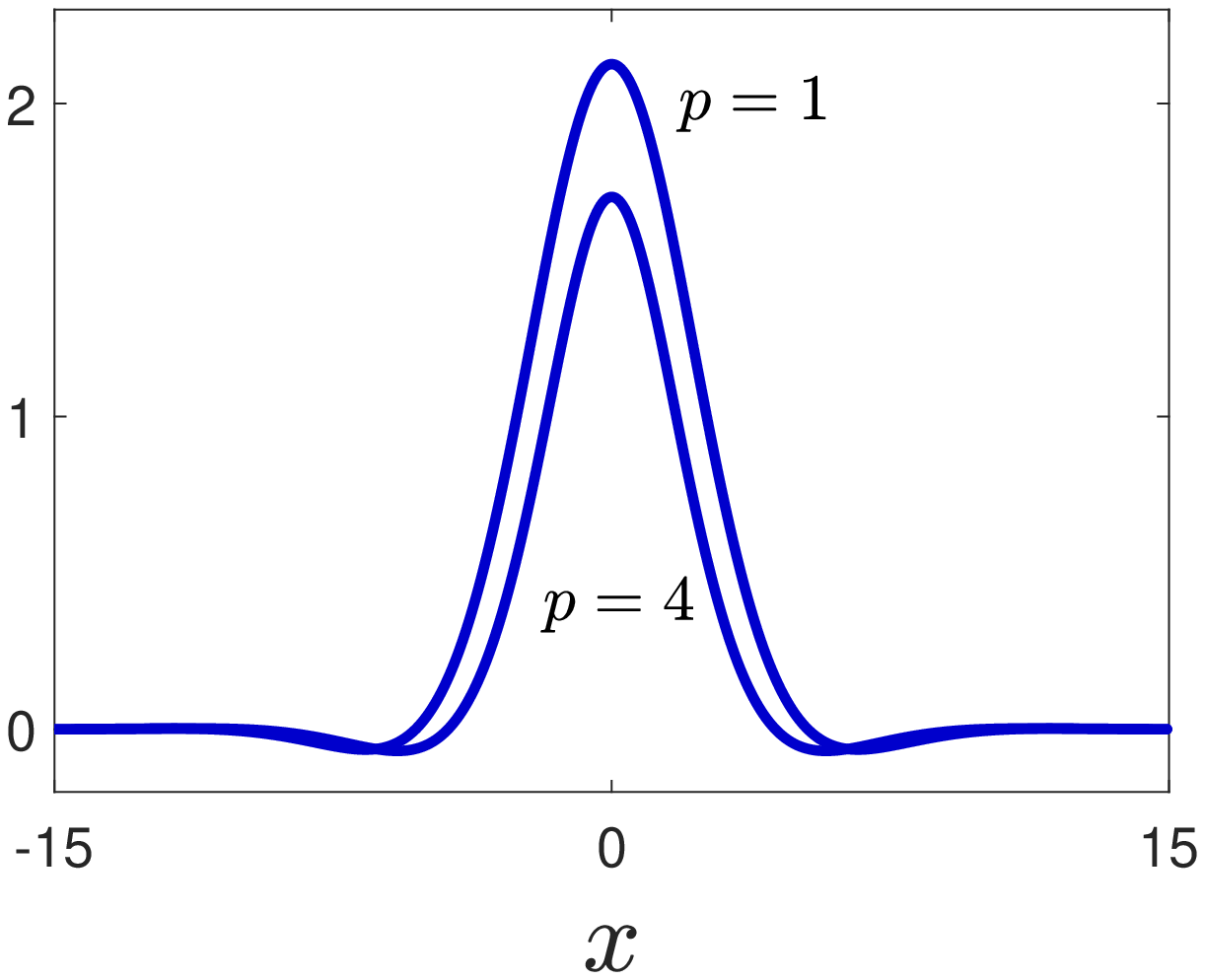}}\hspace*{10pt}
    \subfloat[Residual and iteration errors]{\label{fig:3b}
    \includegraphics[height=0.24\textheight,scale=1.50,keepaspectratio]{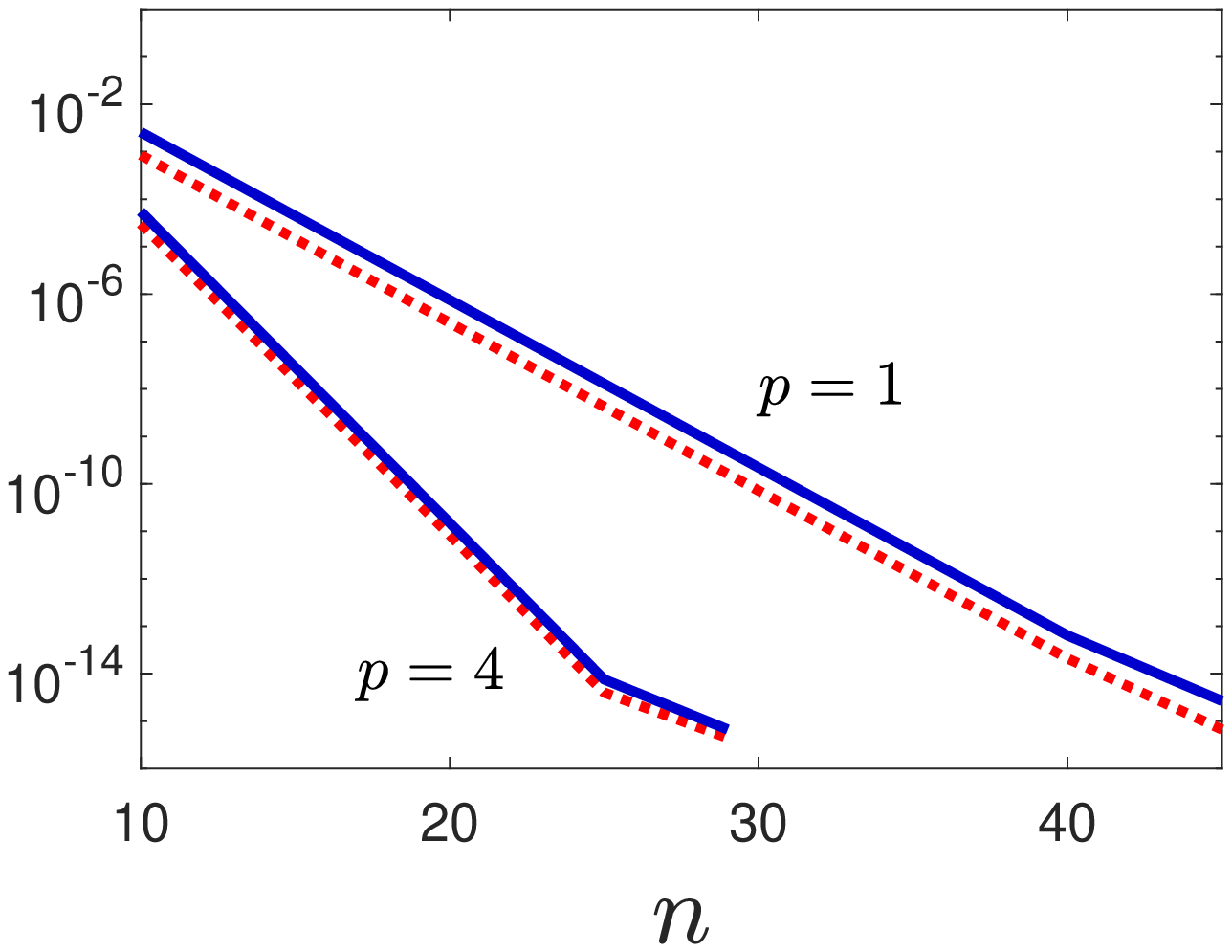}}
    \caption{(a) The numerical solution profiles  of the Rosenau equation (quadratic nonlinearity ($p=1$) and quintic nonlinearity ($p=4$)) for the Gaussian starting function.    (b) Variation of  the iteration  error $E_{n}^{(s)}$ and the residual error $E_{n}^{(r)}$ with the iteration number $n$ for the Rosenau equation with quadratic nonlinearity ($p=1$) and  quintic nonlinearity ($p=4$) for the Gaussian starting function. The dotted (red) line and the solid (blue) line represent $E_{n}^{(s)}$ and $E_{n}^{(r)}$, respectively. (The wave speed $c=1.8$, the computational domain $[-15, \, 15]$ and the mesh size   $0.05$ are used in all computations.) }
    \label{fig:Fig3}
\end{figure}

 We now apply the Petviashvili method based on (\ref{eq:Pstab}) and (\ref{eq:meth2})  to get solitary waves of the Rosenau equation (\ref{eq:re}) with $g(u)=u^{p+1}/(p+1)$. In all the numerical experiments related to the Rosenau equation, we fix the wave speed $c$, the size of the computational domain and the mesh size as $c=1.8$, $-15\leq x \leq 15$ and $h=0.05$ (which corresponds to 600 discrete Fourier modes), respectively. We take the initial guess $\psi_{0}$  as the Gaussian function $\displaystyle    \psi_{0}(x)=e^{-x^{2}}$.   In Figure \ref{fig:3a} we plot the numerical solution for both $p=1$ and $p=4$.  It is interesting to note that all the solution profiles in  Figure \ref{fig:3a}  are symmetric but the tails are non-monotonic and that they assume negative values for some range of $x$. This is exactly contrary to the case for the BBM equation, where the $\mbox{sech}$-type solitary wave solution given explicitly by (\ref{eq:solitary}) is a symmetric localized function with monotonic tails and it takes on only nonnegative values (see also the solution profiles in Figures \ref{fig:1a} and  \ref{fig:2a}). This is due to the different characteristics of the linear dispersion relation of the Rosenau equation, compared to that of the BBM equation. Recall that we are looking for solutions with $u \rightarrow 0$ for $x \rightarrow \pm \infty$. Since $u^{p} \ll u$ ($p \ge 2$)  for large $x$, the monotonicity behavior of the tails of the localized wave will be determined by the linearized version of the wave equation.   The linearized version of (\ref{eq:sw}),  $\mathcal{L}\phi= (cD^4+(c-1))\phi= 0$, admits  solutions of the form $\phi(x)=e^{kx}$ if $ck^4+c-1=0$. The four roots of this quartic equation for $k$ are $c_{0}(1\pm i)/\sqrt{2}$ and $c_{0}(-1\pm i)/\sqrt{2}$ with $c_{0}=\left((c-1)/c\right)^{1/4}>0$. The presence of imaginary roots explains why the non-monotonic (oscillatory) behaviors in Figure \ref{fig:3a} appear. For the linearized version of the BBM equation, the traveling wave solutions satisfy $\mathcal{L}\phi= (cD^2-(c-1))\phi= 0$. It admits solutions of the form $\phi(x)=e^{kx}$ if $k^2=(c-1)/c>0$. Since all the roots of this quadratic equation are real, we have the monotonic solution profiles in Figures \ref{fig:1a} and \ref{fig:2a} for the BBM equation.

\begin{figure}[h!]
    \centering
    \subfloat[Amplitude with $p$]{\label{fig:4a}
    \includegraphics[height=0.24\textheight,scale=1.50,keepaspectratio]{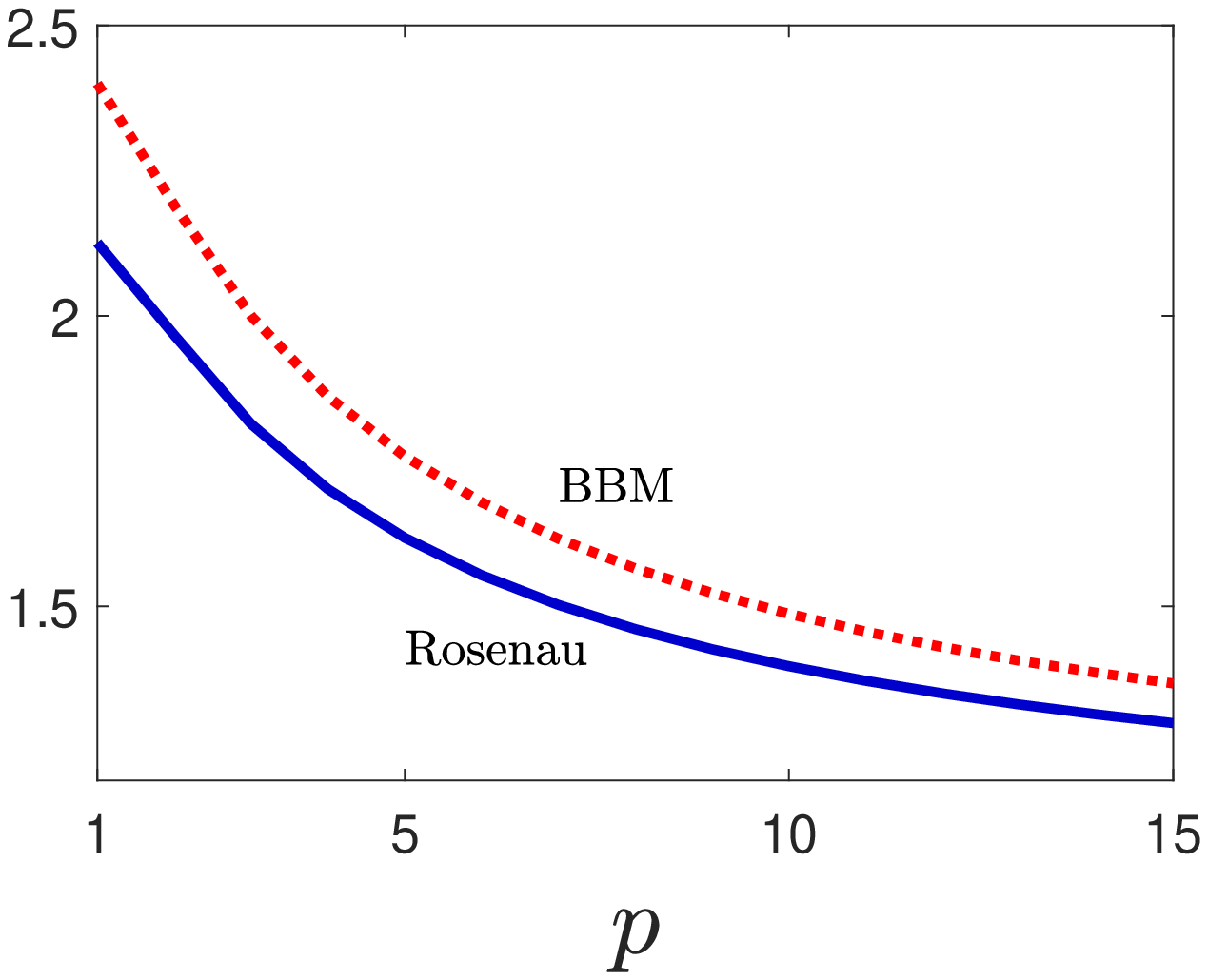}}\hspace*{10pt}
    \subfloat[Amplitude with $\gamma$]{\label{fig:4b}
    \includegraphics[height=0.24\textheight,scale=1.50,keepaspectratio]{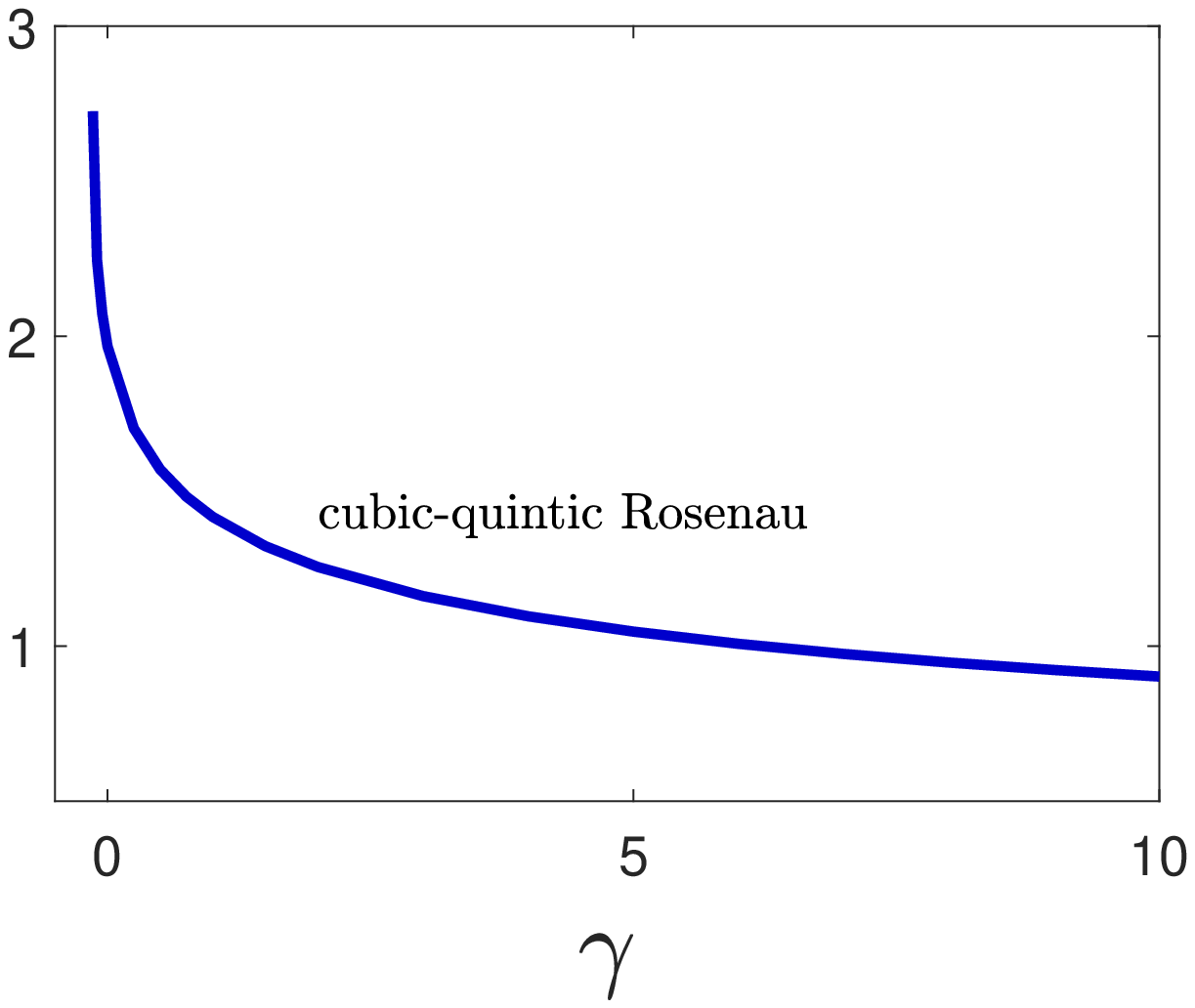}}
    \caption{(a) Variation of the amplitude  of the solitary wave with the degree of nonlinearity ($p$) for the BBM and Rosenau equations with single power nonlinearity.  The dotted (red) line and the solid (blue) line show the exact result for the BBM equation and the numerical result for the Rosenau equation, respectively.   (b) Variation of the amplitude  of the solitary wave with the parameter $\gamma$ for the cubic-quintic Rosenau equation.   (The wave speed $c=1.8$, the computational domain $[-15, \, 15]$ and the mesh size   $0.05$ are used in all computations.) }
    \label{fig:Fig4}
\end{figure}

 In Figure \ref{fig:3b}, using a semi-logarithmic scale, we  present variation of the residual error $E_{r}^{n}$ and the iteration error $E_{i}^{n}$ with the number of iterations ($n$) in the fixed point algorithm. As in the case of the BBM equation, the iteration scheme rapidly converges to the solitary wave solution.    When we make similar calculations for higher values of $p$,  we always get a similar profile but  the amplitude ($\phi(0)$) decreases as  $p$ increases. Figure \ref{fig:4a} shows the variation of the amplitude of the solitary wave with the degree of the nonlinearity, $p$.  In Figure \ref{fig:4a}, for comparison purposes, we also plot the analytical relation between the amplitude $A$ of the $sech$-type solitary wave solution of the BBM equation with $p$. The figure shows that the overall behaviors are very similar for the BBM and the Rosenau equations in the sense that in both cases the amplitudes exhibit monotonic decreasing behavior for increasing values of  $p$, with a horizontal asymptote at 1. We remind the reader that the amplitudes of the solitary wave solutions depend on the chosen value of the wave speed $c$ and that, in our experiments, they are  larger than one  for $c=1.8$. Those observations may not be valid for some other values of $c$.

\begin{figure}[h!]
    \centering
    \subfloat[Numerical profiles]{\label{fig:5a}
    \includegraphics[height=0.24\textheight,scale=1.50,keepaspectratio]{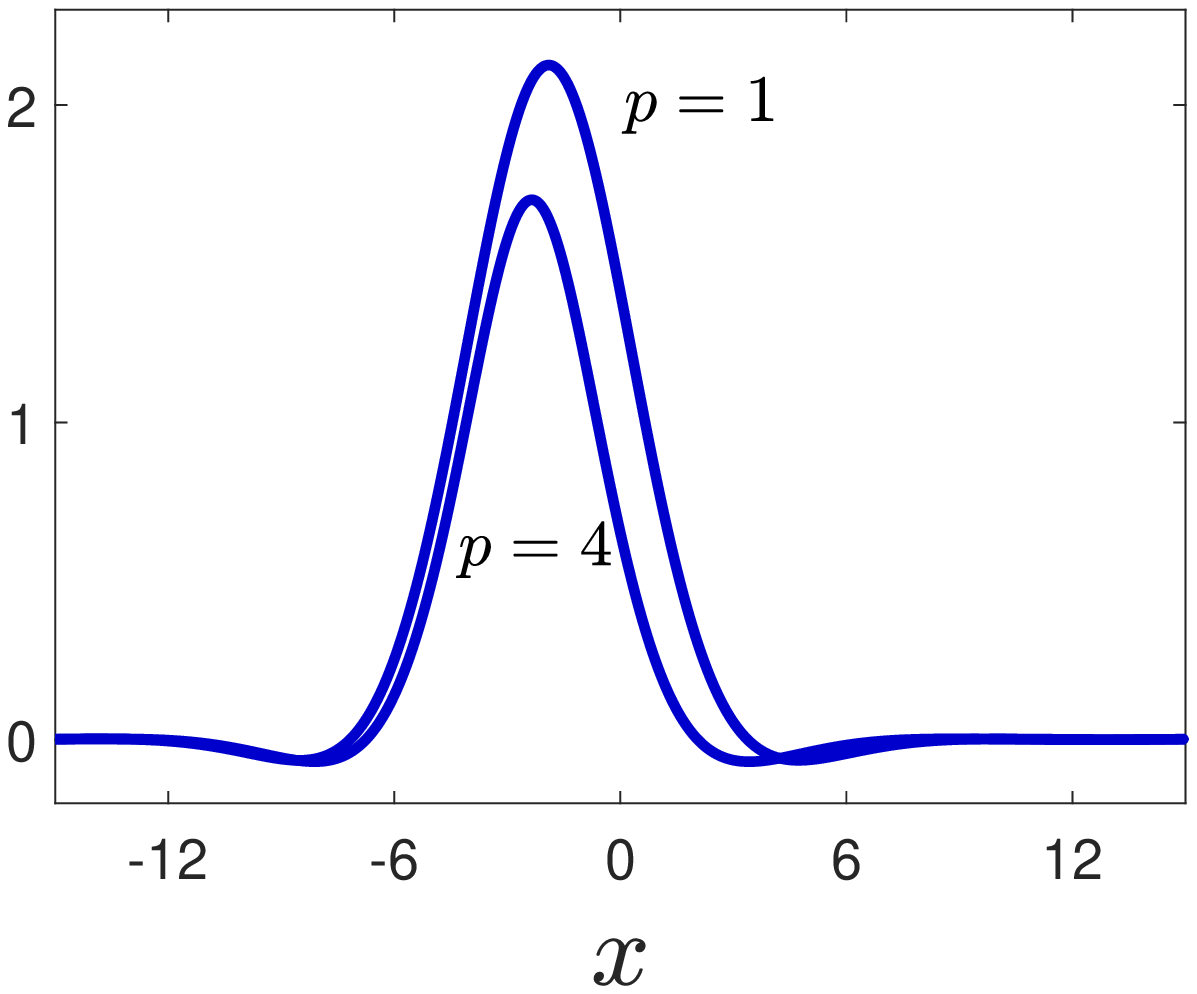}}\hspace*{10pt}
    \subfloat[Residual and iteration errors]{\label{fig:5b}
    \includegraphics[height=0.24\textheight,scale=1.50,keepaspectratio]{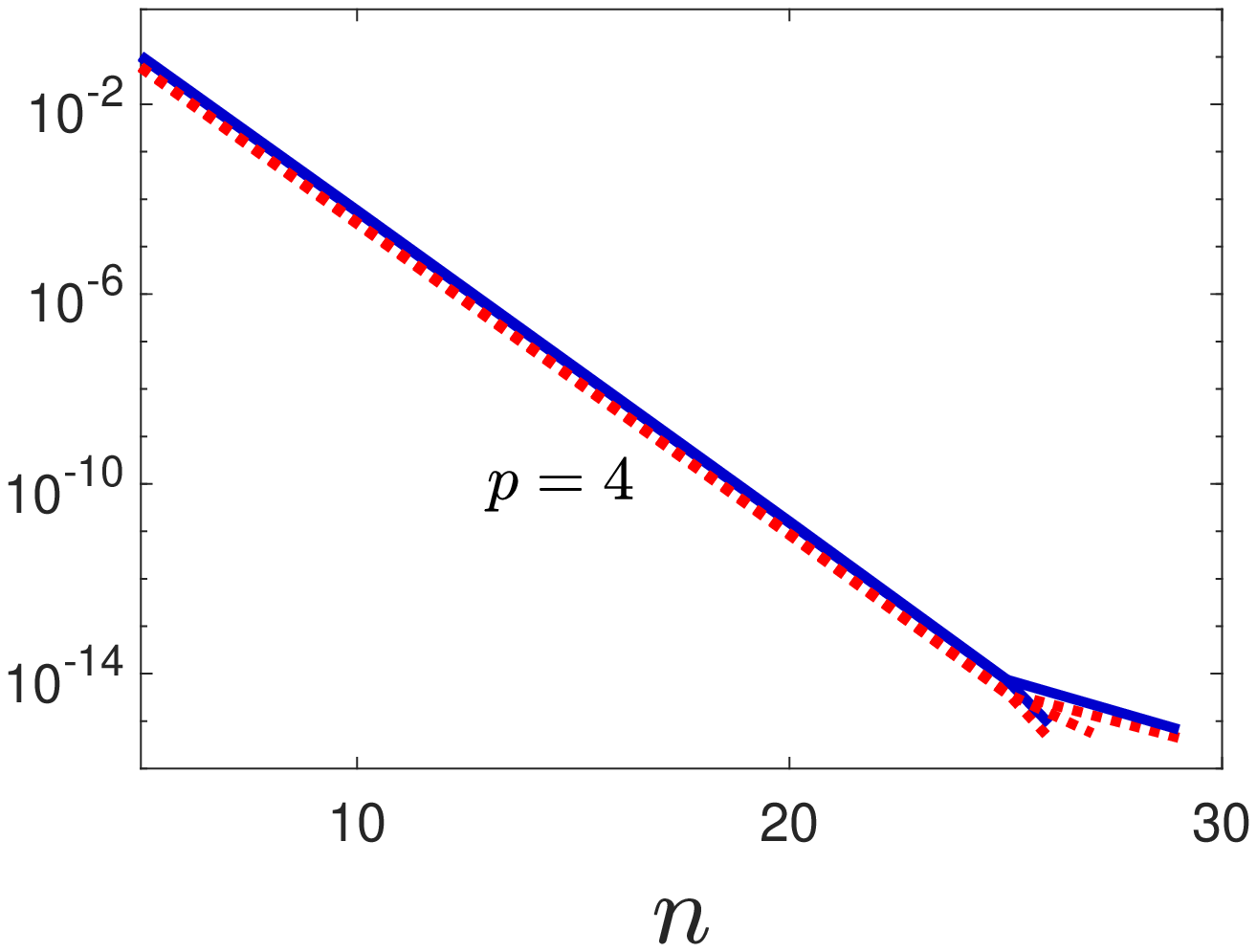}}
    \caption{(a) The numerical solution profiles of the Rosenau equation (quadratic nonlinearity ($p=1$) and quintic nonlinearity ($p=4$)) for the triangular starting function.     (b) Variation of  the iteration  error $E_{n}^{(s)}$ and the residual error $E_{n}^{(r)}$ with the iteration number $n$ for  three different types of the starting function. The Rosenau equation is considered with  quintic nonlinearity ($p=4$). The dotted (red) line and the solid (blue) line represent $E_{n}^{(s)}$ and $E_{n}^{(r)}$, respectively.  (The wave speed $c=1.8$, the computational domain $[-15, \, 15]$ and the mesh size   $0.05$ are used in all computations.) }
    \label{fig:Fig5}
\end{figure}

As in the case of the BBM equation, we now discuss the robustness of the iteration scheme to the initial guess for the case of the Rosenau equation. When the initial guess $\psi_{0}(x)$  is taken as  the Gaussian function which is smooth and symmetric function about $x=0$,  we observe that  the algorithm converges rapidly to the solitary wave solution which  is symmetric  about the origin. When the initial guess is taken as the nonsmooth symmetric function  $\displaystyle  \psi_{0}(x)=e^{-\vert x \vert}$, again we obtain almost the same numerical results and we do not show the corresponding figures here. As the initial guess we now consider the asymmetric triangular function $\psi_{0}(x)$ defined by $\displaystyle  \psi_{0}(x)=x/3+2$ for $-6 \leq x \leq -3$, $\displaystyle  \psi_{0}(x)=-x/9+2/3$ for $-3 < x \leq 6$ and  $\displaystyle  \psi_{0}(x)=0$ otherwise. We get similar results with those obtained for the BBM equation when the starting function is asymmetric triangular function and we plot them in Figure \ref{fig:5a}.  Again, the iteration scheme converges rapidly and the solitary wave solution is a shift of the solitary wave obtained in the previous case. We now compare the algorithmic errors and the residual errors for three different starting functions: Gaussian function, symmetric nonsmooth exponential function, asymmetric triangular function. We present the results for  the residual error $E_{r}^{n}$ and the iteration error $E_{r}^{n}$ with the number of iterations ($n$)  in Figure \ref{fig:5b}. We observe that the curves corresponding to three different starting functions are almost  indistinguishable from each other. We conclude that these experiments illustrate the robustness of the Petviashvili method.

\subsection{The  Cubic-Quintic Rosenau Equation}
 We now apply the Petviashvili method to the Rosenau equation with double power nonlinearity. For simplicity, we take  $g(u)=u^{3}/3+\gamma u^{5}/5$ in (\ref{eq:re}), where $\gamma$ is a constant parameter, and we construct the solitary wave solution numerically for the cubic-quintic Rosenau equation
 \begin{equation}
    u_{t}+u_{x}+u_{xxxxt}+u^2 u_x+\gamma u^4 u_x=0.  \label{eq:double}
\end{equation}
 Traveling wave solutions $u(x,t)=\phi(x-ct)$ of  (\ref{eq:double})  satisfy
\begin{equation}
    \mathcal{L}\phi=cD^4+(c-1))\phi={1\over 3}\phi^{3}+{\gamma \over 5} \phi^{5}  \label{eq:doublesw}
\end{equation}
under the assumption that $\phi(x)$ and its derivatives tend to zero as  $\vert x\vert\to \infty$.
In the Fourier space, (\ref{eq:doublesw}) becomes
\begin{equation}
   l(\xi) \widehat{\phi}(\xi)= {1\over 3}\widehat{\phi^{3}}(\xi)+{\gamma\over 5}\widehat{\phi^{5}}(\xi) , ~~~~l(\xi)=c \xi^4+c-1 \label{eq:doublesw2}
\end{equation}
 where $l(\xi)$ is the symbol of $ \mathcal{L}$.  It is worth pointing out that under the assumption  $c>1$ we have $l(\xi)>0$  for any $\xi \in \mathbb{R}$.  Multiplying (\ref{eq:doublesw2}) by $(\widehat{\phi})^{*}$  and integrating over $\mathbb{R}$ gives
\begin{equation}
    \big\langle  l(\xi)\widehat{\phi}(\xi),(\widehat{\phi})^{*}(\xi)\big\rangle
        = {1\over 3}  \big\langle  \widehat{\phi^{3}}(\xi),  (\widehat{\phi})^{*}(\xi)\big\rangle
        +{\gamma\over 5} \big\langle  \widehat{\phi^{5}}(\xi), (\widehat{\phi})^{*}(\xi)\big\rangle   \label{eq:doublesw3}.
\end{equation}
Assuming that  $\psi_n(x)$ represents the approximation at the $n$th iteration to $\phi(x)$, we now introduce the two stabilizing factors
 \begin{equation*}
    P_n=3{\big\langle  l(\xi)\widehat{\psi_n}(\xi),~(\widehat{\psi_n})^{*}(\xi)\big\rangle \over
        \big\langle  \widehat{\psi_n^{3}}(\xi),~  (\widehat{\psi_n})^{*}(\xi)\big\rangle}, ~~~~~
    Q_n={5\over\gamma}  {\big\langle  l(\xi)\widehat{\psi_n}(\xi),~(\widehat{\psi_n})^{*}(\xi)\big\rangle \over
        \big\langle  \widehat{\psi_n^{5}}(\xi),~  (\widehat{\psi_n})^{*}(\xi)\big\rangle}    .
        \label{eq:PQstab}
\end{equation*}

Because of (\ref{eq:doublesw3}), we note that  $1/P_{n}+1/Q_{n} \rightarrow 1$ as $n \rightarrow \infty$ when the scheme converges, that is, when $\psi_{n} \rightarrow \phi$ as  $n \rightarrow \infty$. To solve (\ref{eq:doublesw2}) numerically we suggest the following iteration scheme
\begin{equation}
    \widehat{\psi_{n+1}}(\xi)= {1\over 3}(P_n)^{3/2} {\widehat{\psi_n^{3}}(\xi) \over l(\xi)}
                                +{\gamma\over 5}(Q_n)^{5/4} {\widehat{\psi_n^{5}}(\xi) \over l(\xi)},   \quad n=0,1,\cdots, \label{eq:doublemeth2}
\end{equation}
where the powers $3/2$ and $5/4$ of $P_{n}$ and $Q_{n}$ are the optimum values corresponding to single power nonlinearity for the cases $p=2$ and $p=4$, respectively. The following numerical experiments show that the sequence $\widehat{\psi}_n(\xi)$ converges to the fixed point of (\ref{eq:doublesw2}) and consequently it converges to the solitary wave solution of the cubic-quintic Rosenau equation (\ref{eq:double}).

\begin{figure}[h!]
    \centering
    \subfloat[Numerical profiles]{\label{fig:6a}
    \includegraphics[height=0.24\textheight,scale=1.50,keepaspectratio]{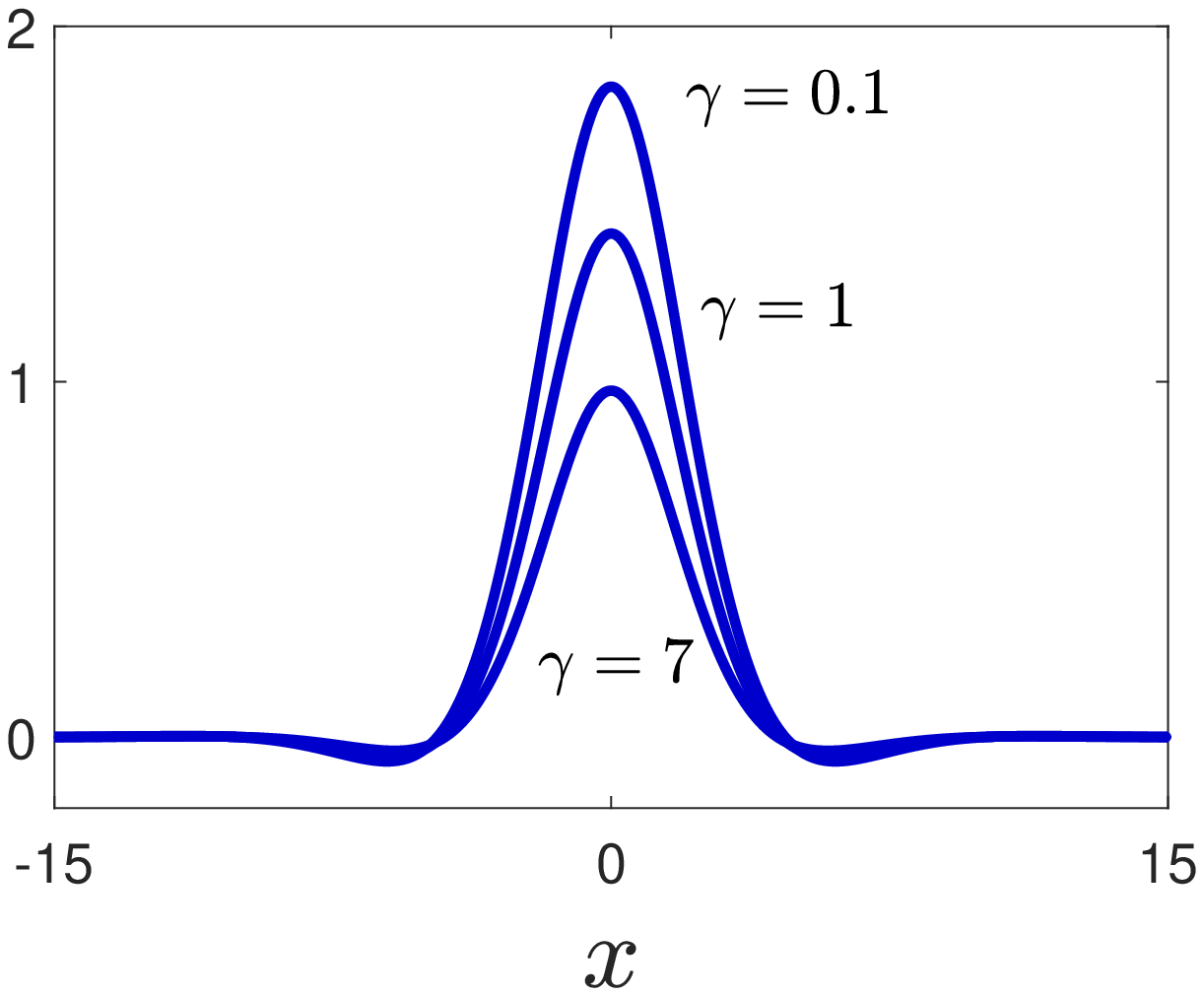}}\hspace*{10pt}
    \subfloat[Residual and iteration errors]{\label{fig:6b}
    \includegraphics[height=0.24\textheight,scale=1.50,keepaspectratio]{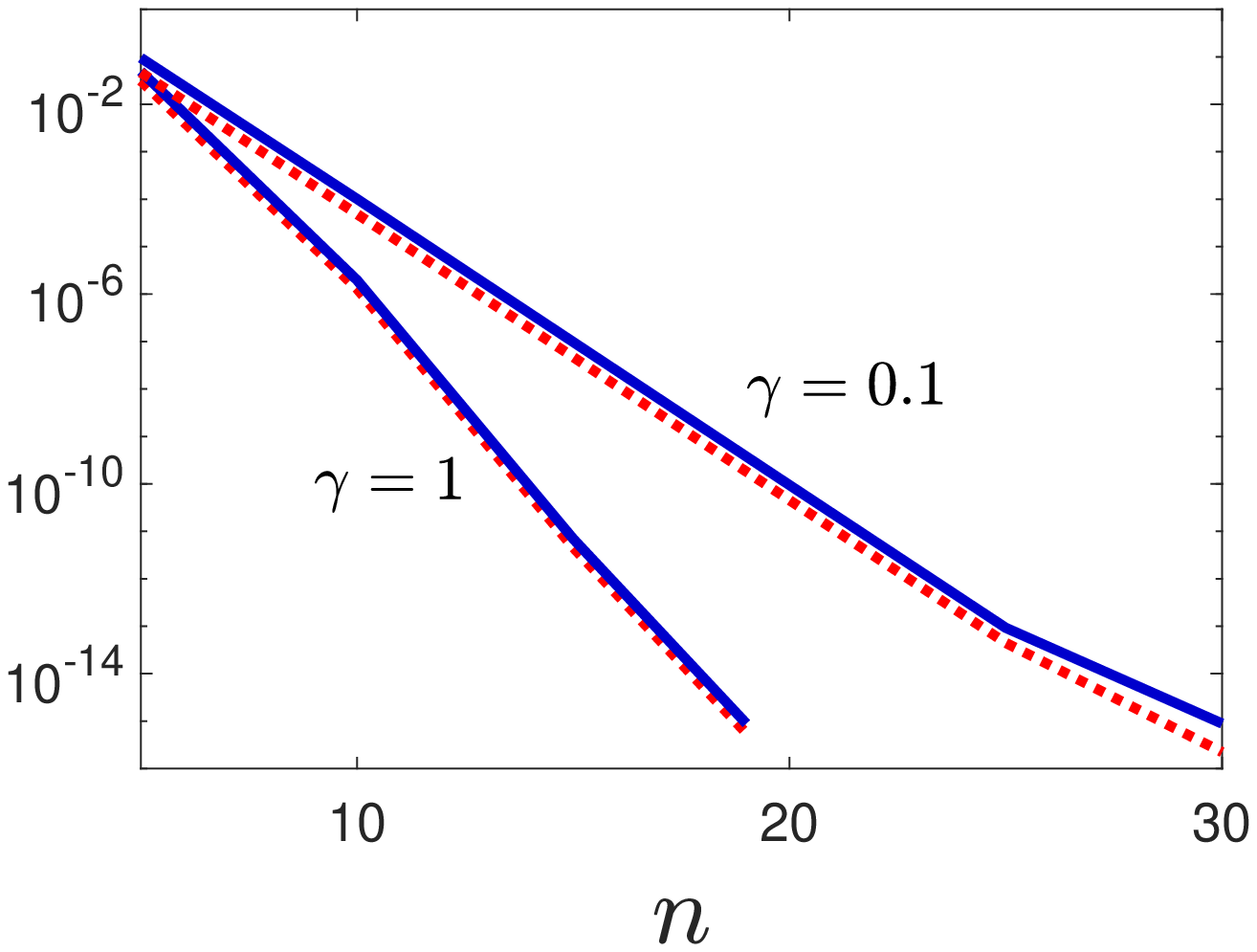}}
    \caption{(a) The numerical solution profiles of the cubic-quintic Rosenau equation  for the Gaussian starting function  and the parameter values $\gamma=0.1$, $\gamma=1$ and $\gamma=7$.     (b) Variation of  the iteration  error $E_{n}^{(s)}$ and the residual error $E_{n}^{(r)}$ with the iteration number $n$ for  $\gamma=1$ and $\gamma=0.1$. The dotted (red) line and the solid (blue) line represent $E_{n}^{(s)}$ and $E_{n}^{(r)}$, respectively. (The wave speed $c=1.8$, the computational domain $[-15, \, 15]$ and the mesh size   $0.1$ are used in all computations.)}
    \label{fig:Fig6}
\end{figure}

Again we take the initial guess $\psi_{0}$  as the Gaussian function $\displaystyle    \psi_{0}(x)=e^{-x^{2}}$ and set $c=1.8$.    In Figure \ref{fig:6a} we plot the profiles of the numerical solutions for  $\gamma=0.1$, $\gamma=1$ and $\gamma=7$. As in the case of single power nonlinearity, all the solution profiles in  Figure \ref{fig:6a}  are symmetric but the tails are non-monotonic and that they assume negative values for some range of $x$. For the present experiment we point out that the   residual error $E_{n}^{(r)}$ and  the iteration error $E_{n}^{(i)}$ based on the stabilizing factor  are defined as
\begin{equation}
    E_{n}^{(r)}= \Big \Vert \mathcal{L}\psi_{n}(x)-\big({1\over 3}\psi_{n}^{3}(x)+{\gamma \over 5} \psi_{n}^{5}(x)\big) \Big \Vert_{L^{\infty}}, ~~~~
    E_{n}^{(s)}= \Big \vert 1-\big({1\over P_{n}}+{1\over Q_{n}}\big)\Big \vert . \label{eq:doubleerror}
\end{equation}
 In Figure \ref{fig:6b}, to simplify the figure we present the variation of errors for $\gamma=1$ and $\gamma=0.1$ only. Again, a semi-logarithmic scale is used in the figure and the variation of the residual error $E_{n}^{(r)}$ and the iteration error $E_{n}^{(s)}$ with the number of iterations ($n$) in the fixed point algorithm is presented.  As in the previous cases, we observe very fast convergence of the iteration.

 We remark that, even though both of (\ref{eq:re-q}) and (\ref{eq:double})  have both cubic and quintic nonlinearities with different coefficients ,  the above-mentioned behaviors of the solution profiles in Figure \ref{fig:6a} (that is, nonmonotonic behavior and taking on negative values) are in sharp contrast to that of the solitary wave solution (\ref{eq:sech}) of  (\ref{eq:re-q}). At this point, we would like to remind the reader that $c=1/2$ in (\ref{eq:sech}) and that our numerical computations are based on the assumption $c>1$.

In Figure \ref{fig:4b} we present the variation of the amplitude ($\phi(0)$) with the parameter $\gamma$ for the traveling wave solutions of (\ref{eq:double}). We observe that the amplitude  decreases as  $\gamma$ increases. We also observe numerically that there is a threshold $\gamma=-0.1395$ below which the Petviashvili method does not converge.

\bibliographystyle{elsarticle-num}
\bibliography{RosenauRef}

\end{document}